








\documentstyle{amsppt}
\magnification=\magstep1
\NoRunningHeads
\pageno=1

\vsize=7.4in


\def\tr{\text{tr }}

\def\dim{\text{dim }}

\def\rank{\text{rank }}


\topmatter
\title
A note on pairs of projections
\endtitle

\author
N.J. Kalton \\
University of Missouri-Columbia
\endauthor

\address
Department of Mathematics,
University of Missouri
Columbia, MO  65211, U.S.A.
\endaddress
\email nigel\@math.missouri.edu \endemail
\thanks Supported by NSF grant DMS-9500125
 \endthanks
\abstract

We give a brief proof of a recent result of Avron, Seiler and Simon.

\endabstract

\endtopmatter

\document \baselineskip=14pt

In \cite{1}, it is proved that if $P,Q$ are (not necessarily
self-adjoint) projections on a Hilbert space and $(P-Q)^n$
is trace-class (i.e. nuclear) for some odd integer $n$ then
$\tr(P-Q)^n$ is an integer and in fact, if $P$ and $Q$ are
self-adjoint, $ \tr(P-Q)^n=\dim E_{10}-\dim E_{01} $ where
$E_{ab}=\{x:Px=ax,\ Qx=bx\};$ (see also \cite{2}).  The
proof given in \cite{1} uses the structure of the spectrum
of $P-Q$ and Lidskii's theorem; it is therefore not
applicable to more general Banach spaces.
The purpose of this note is to give a very brief proof of the same result
which involves only simple algebraic identities and is valid
in any Banach space with a well-defined trace (i.e. with the
approximation property).  We use $[A,B]$ to denote the
commutator $AB-BA.$

The basic material about operators on Banach spaces which we use can be
found in the book of Pietsch \cite{3}.  We summarize the two most
important ingredients.

We  will  need  the  following  basic  result  from Fredholm
theory.  Suppose $X$ is a Banach space and $A:X\to X$ is  an
operator such  that for  some $m,$  $A^m$ is  compact.   Let
$S=I-A$;     then     $F=\cup_{k\ge     1}S^{-k}(0)$      is
finite-dimensional and if  $Y=\cap_{k\ge 1}S^k(X)$ then  $Y$
is  closed  and  $X$  can  be  decomposed  as  a  direct sum
$X=F\oplus Y.$ Furthermore $F$ and $Y$ are invariant for $S$
and $S$ is invertible on $Y.$ We refer to \cite{3} 3.2.9 (p.
141-142) for a slightly more general result.

We  will  also  need  the  following  properties  of nuclear
operators and the trace.  If  $X$ is a Banach space then  an
operator $T:X\to X$ is called nuclear if it can written as a
series $T=\sum_{n=1}^{\infty}A_n$ where each $A_n$ has  rank
one  and  $\sum_{n=1}^{\infty}\|A_n\|<\infty.$  The  nuclear
operators form an ideal  in the space of  bounded operators.
When $X$ has the approximation property, one can then define
the   trace   of    $T$   {\it   unambiguously}    by   $\tr
T=\sum_{n=1}^{\infty}\tr A_n$ (where the trace of a rank one
operator $A=x^*\otimes  x$ is  defined in  the usual  way by
$\tr A=x^*(x).$) The  trace is then  a linear functional  on
the ideal  of nuclear  operators and  has the  property that
$\tr[A,T]=0$ if  $A$ is  bounded and  $T$ is  nuclear.   See
Chapter 4 of \cite{3} and particularly Theorem 4.7.2.

\proclaim{Lemma} Let $X$ be  a Banach space and  suppose $P$
and  $Q$  are  two   projections  on  $X$.     Let  $M=P-Q$,
$U=(I-Q)(I-P)+QP$, $V=(I-P)(I-Q)+PQ$ and suppose $T$ is  any
operator which commutes with both $P$ and $Q$.  Then\newline
(1)  $M^2$  commutes  with  both  $P$  and  $Q$.\newline (2)
$[(I-2Q)TM,PV]=TM(I-M^2).$\newline   (3)   If   $I-M^2$   is
invertible $[(I-2Q)TM(I-M^2)^{-1},PV]=TM.$\endproclaim

\demo{Proof}(1) was first  observed by Dixmier,  Kadison and
Mackey  as  remarked  in  \cite{1}.    For  (2) observe that
$QU=UP$        and        $UV=VU=        I-M^2.$       Hence
$M(I-M^2)=PUV-QUV=PVU-UPV=-[U,PV]=[I-U,PV]=[(I-2Q)M,PV].$ If
$T$ commutes with $P$ and  $Q$ then (2) follows.   Note that
(3)   is    immediate   from    (2),   replacing    $T$   by
$T(I-M^2)^{-1}.$\qed\enddemo

\proclaim{Theorem}Let  $X$  be  a  Banach  space  with   the
approximation property, and suppose  $n$ is an odd  integer.
If $P,Q$ are two projections  on $X$ so that $(P-Q)^{n}$  is
nuclear,     then     $\tr(P-Q)^n=\dim     E_{10}-\dim\tilde
E_{01}=\dim\tilde E_{10}-\dim E_{01},$ where  $E_{ab}=\{x\in
X:Px=ax,\    Qx=bx\}$     and    $\tilde     E_{ab}=\{x^*\in
X^*:P^*x^*=ax^*,\ Q^*x^*=bx^*\}.$\endproclaim

\demo{Remark}If $X$ is a Hilbert  space and $P$ and $Q$  are
self-adjoint  this  is  equivalent  to  the result of Avron,
Seiler and Simon.\enddemo

\demo{Proof}We use the notation of  the lemma.  If $M^n$  is
nuclear then some power of $M^2$ is compact.  Let  $S=I-M^2$
and  let  $F=\cup_{k\ge  1}  S^{-k}(0)$  and   $Y=\cap_{k\ge
1}S^k(X).$ Then as noted above we have that $\dim F<\infty$,
$X=F\oplus Y$ and  $S$ is invertible  on $Y.$ Since  $P$ and
$Q$ commute with $S$ both $F$ and $Y$ are invariant for  $P$
and $Q$.   We denote the  restriction of an  operator $T$ to
$F$ or $Y$ by $T_F$ or $T_Y.$

Now $(I-M^2_Y)$  is invertible  on $Y$  so that  (3) of  the
lemma  expresses  $M_Y^n$  as  the  commutator  of a nuclear
operator and a bounded operator.  Hence $\tr M_Y^n=0.$

On the  other hand,  by (2)  of the  Lemma, $M_F-M_F^n$ is a
commutator on $F$ which  is finite-dimensional so that  $\tr
M_F^n=\tr M_F=\tr P_F-\tr Q_F \in \Bbb Z.$

It is  easy to  see from  elementary computations  that $\tr
P_F-\tr Q_F=  \dim F-\rank  (I-P_F)-\rank Q_F  =\dim F- \dim
((I-P)F+Q(F)) -\dim E_{01}.$ Now $\dim  F-\dim((I-P)F+Q(F))$
is the dimension of the subspace of $F^*$ of all $f^*$  such
that $P_F^*f^*=f^*$ and  $Q_F^*f^*=0;$ if we  identify $F^*$
with $Y^{\perp}$ via the direct sum decomposition this space
coincides with $\tilde E_{10}.$ Now, it follows easily  from
the  properties  of  the  trace  that $\tr M^n=\tr M_F^n+\tr
M_Y^n$.  This  gives the second  formula for $\tr  M^n.$ The
other formula is similar.  \qed\enddemo

\enddocument

\bye